\documentclass[a4paper, 11pt, twoside]{article}
\usepackage{amsmath}
\usepackage{amscd}
\usepackage{amssymb}
\usepackage{amsthm}
\usepackage{graphics}
\usepackage{indentfirst}
\usepackage{latexsym}
\usepackage{amsfonts}
\usepackage{multicol}
\usepackage{stmaryrd}
\usepackage{array}
\usepackage{pxfonts}
\usepackage{verbatim}
\usepackage{bbm}

\makeatletter
\newenvironment{pf}[1][\proofname] {\par\pushQED{\qed}\normalfont\topsep6\p@\@plus6\p@\relax\trivlist\item[\hskip\labelsep\bfseries#1\@addpunct{.}]\ignorespaces}{\popQED\endtrivlist\@endpefalse}
\makeatother

\newtheoremstyle{mattthm}{}{}{\itshape}{}{\bfseries}{.}{ }{}

\theoremstyle{mattthm}
\newtheorem{lemma}{Lemma}[section]
\newtheorem{propn}[lemma]{Proposition}
\newtheorem{thm}[lemma]{Theorem}
\newtheorem{cory}[lemma]{Corollary}

\newtheoremstyle{mattdef}{}{}{}{}{\bfseries}{.}{ }{}

\theoremstyle{mattdef}

\newtheorem*{eg}{Example}

\begin{document}

\hyphenation{multi-reg-ular}
\newenvironment{pfenum}{\begin{pf}\indent\begin{enumerate}\vspace{-\topsep}}{\qedhere\end{enumerate}\end{pf}}

\newcommand\reg{\mathcal{R}}
\newcommand\swed[1]{\stackrel{#1}{\wed}}
\newcommand\inv[2]{\llbracket#1,#2\rrbracket}
\newcommand\one{\mathbbm{1}}
\newcommand\bsm{\begin{smallmatrix}}
\newcommand\esm{\end{smallmatrix}}
\newcommand{\rt}[1]{\rotatebox{90}{$#1$}}
\newcommand\la\lambda
\newcommand{\ol}{\overline}
\newcommand{\ul}{\underline}
\newcommand{\lan}{\langle}
\newcommand{\ran}{\rangle}
\newcommand\partn{\mathcal{P}}
\newcommand\calp{\mathcal{P}}
\newcommand\calf{\mathcal{F}}
\newcommand{\py}[3]{\,_{#1}{#2}_{#3}}
\newcommand{\pyy}[5]{\,_{#1}{#2}_{#3}{#4}_{#5}}
\newcommand{\thmlc}[3]{\textup{\textbf{(\!\! #1 \cite[#3]{#2})}}}
\newcommand{\sss}{\mathfrak{S}_}
\newcommand{\dom}{\trianglerighteqslant}
\newcommand{\doms}{\vartriangleright}
\newcommand{\ndom}{\ntrianglerighteqslant}
\newcommand{\ndoms}{\not\vartriangleright}
\newcommand{\domby}{\trianglelefteqslant}
\newcommand{\domsby}{\vartriangleleft}
\newcommand{\ndomby}{\ntrianglelefteqslant}
\newcommand{\ndomsby}{\not\vartriangleleft}
\newcommand{\subs}[1]{\subsection{#1}}
\newcommand{\nin}{\notin}
\newcommand{\nchar}{\operatorname{char}}
\newcommand{\thmcite}[2]{\textup{\textbf{\cite[#2]{#1}}}\ }
\newcommand\zez{\mathbb{Z}/e\mathbb{Z}}
\newcommand\zepz{\mathbb{Z}/(e+1)\mathbb{Z}}
\newcommand{\bbf}{\mathbb{F}}
\newcommand{\bbc}{\mathbb{C}}
\newcommand{\bbn}{\mathbb{N}}
\newcommand{\bbq}{\mathbb{Q}}
\newcommand{\bbz}{\mathbb{Z}}
\newcommand\zo{\bbn_0}
\newcommand{\gs}{\geqslant}
\newcommand{\ls}{\leqslant}
\newcommand\dw{^\triangle}
\newcommand\wod{^\triangledown}
\newcommand{\hhh}{\mathcal{H}_}
\newcommand{\sect}[1]{\section{#1}}
\newcommand{\ff}{\mathfrak{f}}
\newcommand{\fff}{\mathfrak{F}}
\newcommand\cf{\mathcal{F}}
\newcommand\fkn{\mathfrak{n}}
\newcommand\sx{x}
\newcommand\bra[1]{|#1\ran}
\newcommand\arb[1]{\widehat{\bra{#1}}}
\newcommand\foc[1]{\mathcal{F}_{#1}}
\newcommand{\clam}{\begin{description}\item[\hspace{\leftmargin}Claim.]}
\newcommand{\prof}{\item[\hspace{\leftmargin}Proof.]}
\newcommand{\malc}{\end{description}}
\newcommand\ppmod[1]{\ (\operatorname{mod}\ #1)}
\newcommand\wed\wedge
\newcommand\wede\barwedge
\newcommand\uu[1]{\,\begin{array}{|@{\,}c@{\,}|}\hline #1\\\hline\end{array}\,}
\newcommand{\ux}[1]{\operatorname{ht}_{#1}}
\newcommand\erim{\operatorname{rim}}
\newcommand\mire{\operatorname{rim}'}
\newcommand\mmod{\ \operatorname{Mod}}
\newcommand\cgs\succcurlyeq
\newcommand\cls\preccurlyeq
\newcommand\cg\succ
\newcommand\cl\prec
\newcommand\inc{\mathfrak{A}}
\newcommand\fsl{\mathfrak{sl}}
\newcommand\ba{\mathbf{s}}
\newcommand\ta{\tilde\ba}
\newcommand\kt[1]{|#1\rangle}
\newcommand\tk[1]{\langle#1|}
\newcommand\ket[1]{s_{#1}}
\newcommand\jn\diamond
\newcommand\UU{\mathcal{U}}
\newcommand\MM[1]{M^{\otimes#1}}
\newcommand\add{\operatorname{add}}
\newcommand\rem{\operatorname{rem}}
\newcommand\La\Lambda
\newcommand\lra\longrightarrow
\newcommand\res{\operatorname{res}}
\newcommand\lexg{>_{\operatorname{lex}}}
\newcommand\lexgs{\gs_{\operatorname{lex}}}
\newcommand\lexl{<_{\operatorname{lex}}}
\newcommand\tru[1]{{#1}_-}
\newcommand\ste[1]{{#1}_+}
\newcommand\out{^{\operatorname{out}}}
\newcommand\lad{\mathcal{L}}
\newcommand\hsl{\widehat{\mathfrak{sl}}}
\newcommand\fkh{\mathfrak{h}}
\newcommand\GG{H}

\title{An LLT-type algorithm for computing\\ higher-level canonical bases}
\author{Matthew Fayers\\\normalsize Queen Mary, University of London, Mile End Road, London E1 4NS, U.K.\\\texttt{\normalsize m.fayers@qmul.ac.uk}}
\date{}
\maketitle
\begin{center}
\end{center}
\markboth{Matthew Fayers}{An LLT-type algorithm for computing higher-level canonical bases}
\pagestyle{myheadings}

\begin{abstract}
We give a fast algorithm for computing the canonical basis of an irreducible highest-weight module for $U_q(\widehat{\mathfrak{sl}}_e)$, generalising the LLT algorithm.
\end{abstract}


\sect{Introduction}

Let $e\gs2$ be an integer.  In this paper we consider the integrable representation theory of the quantised enveloping algebra $\UU=U_q(\widehat{\mathfrak{sl}}_e)$.  For any dominant integral weight $\La$ for $\UU$, the irreducible highest-weight module $V(\La)$ for $\UU$ can be constructed as a submodule $M^\ba$ of a \emph{Fock space} $\calf^\ba$ (which depends not just on $\La$ but on an ordering of the fundamental weights involved in $\La$).  Using the standard basis of the Fock space, one can define a \emph{canonical basis} (in the sense of Lusztig/Kashiwara) for $M^\ba$.  There is considerable interest in computing this canonical basis (that is, computing the transition coefficients from the canonical basis to the standard basis) because of Ariki's theorem, which says that these coefficients, evaluated at $q=1$, yield decomposition numbers for certain cyclotomic Hecke algebras.  In the case where $\La$ is of level $1$, there is a fast algorithm due to Lascoux, Leclerc and Thibon \cite{llt} for computing the canonical basis.  The purpose of this paper is to give a generalisation of this algorithm to higher levels.

Leclerc and Thibon \cite{lt} showed how the canonical basis could be extended to a basis for the whole of the Fock space in the level $1$ case.  This was generalised to higher levels by Uglov, but using a `twisted' Fock space (which is not obviously isomorphic to a tensor product of level $1$ Fock spaces).  By using Uglov's construction and taking a limit, one can define a canonical basis for the whole of the (untwisted) Fock space, and this in principle gives an algorithm for computing the canonical basis of $M^\ba$.  However, in practice this algorithm is extremely slow.  We give a much faster algorithm here; the way we do this is to compute the canonical basis for an intermediate module $M^{\otimes\ba}$, which is defined to be the tensor product of level $1$ highest-weight irreducibles.  It is then straightforward to discard unwanted vectors to get the canonical basis for $M^\ba$.

We remark that Jacon \cite{j} and Yvonne \cite{yv2} have also given algorithms for computing higher-level canonical bases.  However, Yvonne's algorithm is very slow, since it computes the canonical basis for the whole of the Fock space, while Jacon's algorithm works in a particular type of twisted Fock space, whereas our algorithm remains in the more natural setting of the untwisted Fock space; although these Fock spaces are isomorphic, so that in principle one canonical basis determines the other, in practice it is very difficult to give an explicit isomorphism.

In the next section we give some basic combinatorial and algebraic background, and establish notation.  In Section \ref{keysec}, we describe in detail how the bar involution on a twisted Fock space is computed, and prove an important property of the bar involution which lies at the heart of our algorithm.  In Section \ref{algsec} we describe our algorithm, and prove that it works.  In Section \ref{egsec}, we give examples, and make some further remarks; these concern the generalisation to the case $e=\infty$, and a brief discussion of how to pass from the canonical basis for $M^{\otimes\ba}$ to the canonical basis for $M^\ba$.  Appendix \ref{app} consists of an index of notation.

\sect{Background}\label{backsec}

\subsection{Some elementary notation}

Throughout this paper, $e$ denotes an integer greater than or equal to $2$ (except in Section \ref{infsec} where we consider the generalisation to the case $e=\infty$).  We write $I$ to denote the set $\zez$, which is used as the indexing set for the Cartan matrix of $\UU$.

For any integers $a\ls b$, we write $\inv ab$ for the `integer interval' $\{a,a+1,\dots,b\}$.

\subsection{Partitions and multipartitions}

A \emph{partition} is a sequence $\la=(\la_1,\la_2,\dots)$ of non-negative integers such that $\la_1\gs\la_2\gs\dots$ and the sum $|\la|=\la_1+\la_2+\dots$ is finite.  We write $\calp$ for the set of all partitions.  The partition $(0,0,\dots)$ is usually written as $\varnothing$.

Now suppose $r\in\bbn$.  An \emph{$r$-multipartition} is an ordered $r$-tuple $\la=(\la^{(1)},\dots,\la^{(r)})$ of partitions.  We write $\calp^r$ for the set of $r$-multipartitions.  For $\la\in\calp^r$, we write $|\la|$ for the sum $|\la^{(1)}|+\dots+|\la^{(r)}|$.  We write $\varnothing^r$ for the $r$-multipartition $(\varnothing,\dots,\varnothing)$.  We shall abuse notation slightly in this paper by not distinguishing between a partition and a $1$-multipartition.

We impose a partial order (the \emph{dominance order}) on $\calp^r$ by saying that $\la$ dominates $\mu$ (and writing $\la\dom\mu$) if we have
\[\sum_{l=1}^{k-1}|\la^{(l)}|+\sum_{i=1}^j\la^{(k)}_i\gs\sum_{l=1}^{k-1}|\mu^{(l)}|+\sum_{i=1}^j\mu^{(k)}_i\]
for each $k\in\inv1r$ and $j\gs1$.

Throughout this paper, we use the following notation for multipartitions.  If $\la=(\la^{(1)},\dots,\la^{(r)})$ is an $r$-multipartition for $r>1$, then we write $\tru\la$ for the ($r-1$)-multipartition $(\la^{(2)},\dots,\la^{(r)})$.  If $\nu$ is an ($r-1$)-multipartition, we write $\ste\nu$ for the $r$-multipartition $(\varnothing,\nu^{(1)},\dots,\nu^{(r-1)})$.  Finally, if $\mu$ is an $r$-multipartition, we write $\mu_0$ for the $r$-multipartition $\ste{(\tru\mu)}=(\varnothing,\mu^{(2)},\dots,\mu^{(r)})$.

If $\la\in\calp^r$, the \emph{Young diagram} of $\la$ is the set
\[[\la] = \left\{(i,j,k)\in\bbn^2\times\inv 1r\ \left|\ j\ls \la^{(k)}_i\right.\right\}.\]
We refer to elements of the set $\bbn^2\times\inv 1r$ as \emph{nodes}, and elements of $[\la]$ as nodes of $\la$.  A node $\fkn$ of $\la$ is \emph{removable} if $[\la]\setminus\{\fkn\}$ is again the Young diagram of a multipartition (we denote this partition $\la_{\fkn}$), while a node $\fkn$ not in $[\la]$ is an \emph{addable node of $\la$} if $[\la]\cup\{\fkn\}$ is the Young diagram of a multipartition (which we denote $\la^{\fkn}$).  We impose a total order on the set of all addable and removable nodes of a multipartition by saying that $(i,j,k)$ is \emph{above} $(i',j',k')$ (or $(i',j',k')$ is \emph{below} $(i,j,k)$) if either $k<k'$ or ($k=k'$ and $i<i'$).

Given $\ba=(s_1,\dots,s_r)\in I^r$, we define the \emph{residue} of a node $(i,j,k)$ to be $j-i+s_k\in I$; if a node has residue $l\in I$, we may refer to it as an $l$-node.  We say that a partition $\la$ is \emph{$e$-regular} if there is no $i$ such that $\la_i=\la_{i+e-1}>0$, and that a multipartition $\la$ is \emph{$e$-multiregular} if $\la^{(k)}$ is $e$-regular for each $k$.  We write $\reg$ for the set of $e$-regular partitions and $\reg^r$ for the set of all $e$-multiregular $r$-multipartitions, if $e$ is understood.

\subsection{The quantum algebra $U_q(\hsl_e)$ and the Fock space}

In this paper, we let $\UU$ denote the quantised enveloping algebra $U_q(\hsl_e)$.  This is a $\bbq(q)$-algebra with generators $e_i,f_i$ for $i\in I$ and $q^h$ for $h\in P^\vee$, where $P^\vee$ is a free $\bbz$-module with basis $\{h_i\mid i\in I\}\cup\{d\}$.  The relations are well known; for example, see \cite[\S4.1]{llt}.  For any integer $m>0$, we write $f_i^{(m)}$ to denote the quantum divided power $f_i^m/[m]!$.  Let $\left\{\left.\La_i\ \right|\ i\in I\right\}$ be a set of fundamental weights for $\UU$.

There are various choices for a comultiplication which makes $\UU$ into a Hopf algebra (and hence allows us to regard the tensor product of two $\UU$-modules as a $\UU$-module).  We use the comultiplication denoted $\Delta$ in \cite{kash}, which is defined by
\begin{align*}
\Delta:e_i &\longmapsto e_i\otimes q^{-h_i}\ +\ 1\otimes e_i,\\
f_i &\longmapsto f_i\otimes 1\ +\ q^{h_i}\otimes f_i,\\
q^{h}&\longmapsto q^{h}\otimes q^{h}
\end{align*}
for all $i\in I$ and all $h\in P^\vee$.

The $\bbq$-linear ring automorphism $\ol{\phantom{o}}:\UU\to\UU$ defined by
\[\ol{e_i}=e_i,\qquad \ol{f_i}=f_i,\qquad \ol q=q^{-1},\qquad \ol{q^h}=q^{-h}\]
for $i\in I$ and $h\in P^\vee$ is called the \emph{bar involution}.

Now we fix $\ba\in I^r$ for some $r\gs1$, and define the \emph{Fock space} $\calf^{\ba}$ to be the $\bbq(q)$-vector space with a basis $\{\ket\la\mid\la\in\calp^r\}$, which we call the \emph{standard basis}.  This has the structure of a $\UU$-module, which we now describe.

Given $\la\in\calp^r$, let $\add_i(\la)$ denote the set of addable $i$-nodes  of $\la$, and $\rem_i(\la)$ the set of removable $i$-nodes.  For each $\fkn\in\add_i(\la)$, define $N(\la,\fkn)$ to be the number of addable $i$-nodes of $\la$ above $\fkn$ minus the number of removable $i$-nodes of $\la$ above $\fkn$.  Now the action of $f_i$ is given by
\[f_i\ket\la = \sum_{\fkn\in\add_i(\la)}q^{N(\la,\fkn)}\ket{\la^\fkn}.\]
Similarly, for each $\fkn\in\rem_i(\la)$, define $M(\la,\fkn)$ to be the number of removable $i$-nodes of $\la$ below $\fkn$ minus the number of addable $i$-nodes of $\la$ below $\fkn$.  The action of $e_i$ is then given by
\[e_i\ket\la = \sum_{\fkn\in\rem_i(\la)}q^{M(\la,\fkn)}\ket{\la_\fkn}.\]
The action of the Cartan subalgebra is given by the statement that $\ket\la$ is a weight vector of weight
\[\La_{s_1}+\dots+\La_{s_r}-\sum_{i\in I}c_i\alpha_i,\]
where $c_i$ denotes the number of $i$-nodes of $\la$.

The Fock space is of interest because the submodule $M^{\ba}$ generated by $\ket{\varnothing^r}$ is isomorphic to the irreducible highest-weight module $V(\La_{s_1}+\dots+\La_{s_r})$.  This submodule inherits a bar involution from $\UU$: this is defined by $\ol{\ket{\varnothing^r}}=\ket{\varnothing^r}$ and $\ol{um} = \ol u\,\ol m$ for all $u\in\UU$ and $m\in M^{\ba}$.  This bar involution allows one to define a \emph{canonical basis} for $M^{\ba}$; this consists of vectors $G^\ba(\mu)$, for $\mu$ lying in some subset of $\calp^r$ (with our conventions, this is what Brundan and Kleshchev \cite{bk} call the set of \emph{regular multipartitions}).  These canonical basis vectors are characterised by the following properties:
\begin{itemize}
\item
$\ol{G^\ba(\mu)}=G^\ba(\mu)$;
\item
if we write $G^\ba(\mu) = \sum_{\la\in\calp^r}d^\ba_{\la\mu}s_\la$ with $d^\ba_{\la\mu}\in\bbq(q)$, then we have $d^\ba_{\mu\mu}=1$ while $d^\ba_{\la\mu}\in q\bbz[q]$ if $\la\neq\mu$.
\end{itemize}
In fact, more is true: the coefficient $d^\ba_{\la\mu}$ is zero unless $\mu\dom\la$ and $\ket\la$ and $\ket\mu$ are weight vectors of the same weight (i.e.\ $\la$ and $\mu$ have the same number of $i$-nodes, for each $i$; in particular, $|\la|=|\mu|$).  Of course, this means that $G^\ba(\mu)$ is a weight vector.

There is considerable interest in computing the canonical basis elements (i.e.\ computing the transition coefficients $d^\ba_{\la\mu}$), because of Ariki's theorem \cite{ari}, which says that the coefficients $d^\ba_{\la\mu}$ specialised at $q=1$ equal decomposition numbers for appropriate cyclotomic Hecke algebras.  In fact, the coefficients $d^\ba_{\la\mu}$ (with $q$ still indeterminate) can be regarded as graded decomposition numbers, thanks to the recent work of Brundan and Kleshchev \cite{bk}.

It is possible to extend the bar involution on $M^{\ba}$ to the whole of $\calf^{\ba}$, as we shall explain below; this yields a canonical basis for the whole of $\calf^\ba$, indexed by the set of all $r$-multipartitions.  Moreover, there is an algorithm to compute this canonical basis, and therefore to compute the canonical basis for $M^{\ba}$, but in practice this is extremely slow.  Our approach is to compute the canonical basis for a module lying in between $M^{\ba}$ and $\calf^{\ba}$.  The way we have defined $\calf^{\ba}$ and our choice of coproduct on $\UU$ mean that there is an isomorphism
\begin{alignat*}2
\calf^{\ba}&\overset{\sim}\lra\,\,&\calf^{(s_1)}&\otimes\dots\otimes\calf^{(s_r)}\\
\intertext{defined by linear extension of}
\ket\la&\longmapsto&\ket{(\la^{(1)})}&\otimes\dots\otimes\ket{(\la^{(r)})}.
\end{alignat*}
We will henceforth identify $\calf^{\ba}$ and $\calf^{(s_1)}\otimes\dots\otimes\calf^{(s_r)}$ via this isomorphism.  Since each $\calf^{(s_k)}$ contains a submodule $M^{(s_k)}$ isomorphic to $V(\La_{s_k})$, $\calf^{\ba}$ contains a submodule $\MM{\ba}=M^{(s_1)}\otimes\dots\otimes M^{(s_r)}$ isomorphic to $V(\La_{s_1})\otimes\dots\otimes V(\La_{s_r})$.  Our algorithm will compute the canonical basis of $\MM{\ba}$.

\subsection{The LLT algorithm}

In this section we restrict attention to the case $r=1$, and explain the LLT algorithm for computing canonical basis elements $G^{(s_1)}(\mu)$.  (In fact, the superscript ${}^{(s_1)}$ is unnecessary here, because $G^{(s_1)}(\mu)$ is independent of $s_1$; in general, $G^\ba(\mu)$ should be unchanged if a fixed element of $I$ is added to $s_1,\dots,s_r$ simultaneously.)  The LLT algorithm was first described in the paper \cite{llt}, to which we refer the reader for more details and examples.

In this section, we write a node $(i,j,1)$ of a $1$-multipartition (i.e.\ a partition) just as $(i,j)$.  For each $l\in\bbn$, we define the $l$th \emph{ladder} in $\bbn^2$ to be the set
\[\lad_l=\left\{\left.(i,j)\in\bbn^2\ \right|\ i+(e-1)(j-1)=l\right\}.\]
All the nodes in $\lad_l$ have the same residue (namely, $s_1+1-l$), and we define the residue of $\lad_l$ to be this residue.  If $\mu$ is a partition, we define the $l$th ladder $\lad_l(\mu)$ of $\mu$ to be the intersection of $\lad_l$ with the Young diagram of $\mu$.

The canonical basis elements for $M^{(s_1)}$ are indexed by the $e$-regular partitions.  To construct $G^{(s_1)}(\mu)$ when $\mu$ is $e$-regular, we begin by constructing an auxiliary vector $A(\mu)$.  Let $l_1<\dots<l_t$ be the values of $l$ for which $\lad_l(\mu)$ is non-empty.  For each $k$, let $a_k$ denote the number of nodes in $\lad_{l_k}(\mu)$, and let $i_k$ denote the residue of $\lad_{l_k}$.  Then the vector $A(\mu)$ is defined by
\[A(\mu) = f_{i_t}^{(a_t)}\dots f_{i_1}^{(a_t)}\ket\varnothing.\]
$A(\mu)$ is obviously bar-invariant, and a lemma due to James \cite[6.3.54 \& 6.3.55]{jk} implies that when we expand $A(\mu)$ as
\[A(\mu) = \sum_{\nu\in\calp}a_\nu\ket\nu,\]
we have $a_\mu=1$, while $a_\nu=0$ unless $\mu\dom\nu$.  This means that $A(\mu)$ must equal $G^{(s_1)}(\mu)$ plus a $\bbq(q+q^{-1})$-linear combination of canonical basis vectors $G^{(s_1)}(\nu)$ with $\mu\doms\nu$.  Assuming (by induction on the dominance order) that these $G^{(s_1)}(\nu)$ have been computed, it is straightforward to subtract the appropriate multiples of these vectors from $A(\mu)$ to recover $G^{(s_1)}(\mu)$.  Moreover, the fact that the coefficients of the standard basis elements in $A(\mu)$ all lie in $\bbz[q,q^{-1}]$ means that the coefficients of the canonical basis elements in $A(\mu)$ lie in $\bbz[q+q^{-1}]$.  A more precise description of the procedure to strip off these canonical basis elements is given in the algorithm in Section \ref{algsec}.

\subsection{Uglov's twisted Fock spaces}

We now return to an arbitrary level $r$, and explain how to extend the bar involution on $M^\ba$ to $\calf^{\ba}$.  This is also done in \cite{bk}, and involves using Uglov's construction \cite{uglov} of twisted Fock spaces, and then taking a limit via Yvonne's theorem \cite{yv}.

Given $\ba\in I^r$ as above, define a \emph{multicharge for $\ba$} to be an $r$-tuple $\ta = (\tilde s_1,\dots,\tilde s_r)\in\bbz^r$ such that $\tilde s_k+e\bbz=s_k$ for each $k$.  Uglov defines a twisted Fock space $\calf^{\ta}$ for each multicharge.  The way this is done is exactly as for $\calf^{\ba}$ above, except that the ordering on the addable and removable nodes of a multipartition is changed: let us define the \emph{integral residue} $\res_{\bbz}(i,j,k)$ of a node $(i,j,k)$ to be $\tilde s_k+j-i$, and then say that the node $(i,j,k)$ is above $(i',j',k')$ if either $\res_{\bbz}(i,j,k)>\res_{\bbz}(i',j',k')$ or ($\res_{\bbz}(i,j,k)=\res_{\bbz}(i',j',k')$ and $k>k'$).  Now the construction of the twisted Fock space $\calf^{\ta}$ is exactly the same as for the Fock space $\calf^\ba$, except for the change of ordering of nodes.  In the case $r=1$, this makes no difference at all, but for higher levels $\calf^{\ta}$ is different; in particular, there is no longer an obvious isomorphism from $\calf^{\ta}$ to a tensor product of level $1$ Fock spaces.

The highest-weight vector $\ket{\varnothing^r}$ in $\calf^{\ta}$ still generates a submodule isomorphic to $V(\La_{s_1}+\dots+\La_{s_r})$, and there is a bar involution on this submodule.  Uglov defines an extension of this bar involution to the whole of $\calf^{\ta}$; this bar involution is compatible with the action of $\UU$ in the sense that $\ol{um}=\ol u\,\ol m$ for all $u\in\UU$ and $m\in\calf^{\ta}$.  Furthermore, if we write
\[\ol{\ket\mu} = \sum_{\la\in\calp^r}b^{\ta}_{\la\mu}\ket\la,\]
then the coefficients $b^{\ta}_{\la\mu}$ satisfy a unitriangularity property which enables the algorithmic construction of a canonical basis for the whole of $\calf^{\ta}$.  We will describe Uglov's bar involution explicitly in the next section.

It is easily seen that if we fix $\la\in\calp^r$ and choose $\tilde\ba$ so that $\tilde s_k-\tilde s_{k+1}$ is large relative to $|\la|$ for each $k$ (certainly $\tilde s_k-\tilde s_{k+1}>|\la|$ is sufficient), then the orderings on the addable and removable nodes of $\la$ are the same, so the action of $\UU$ on $\ket\la$ is the same in $\calf^{\ba}$ as in $\calf^{\ta}$.  So $\calf^{\ba}$ can be viewed as a limit of twisted Fock spaces.  To define the bar involution on $\calf^{\ba}$, we need the following stability property of the coefficients $b^{\ta}_{\la\mu}$.

\begin{thm}\label{yvonne}\thmcite{yv}{Theorem 5.2}
Take $\mu\in\calp^r$.  Then there is an integer $N$ such that if $\tilde s_k-\tilde s_{k+1}\gs N$ for each $k$, the transition coefficients $b^{\ta}_{\la\mu}$ are independent of $\ta$.
\end{thm}

This theorem allows us to define a bar involution on $\calf^{\ba}$: for any $\mu$, we choose a multicharge $\ta$ such that $\tilde s_k-\tilde s_{k+1}$ is large relative to $\mu$ for each $k$, and set $b^{\ba}_{\la\mu} = b^{\ta}_{\la\mu}$ for each $\la$.  Then we define
\[\ol{\ket\mu} = \sum_{\la\in\calp^r}b^{\ba}_{\la\mu}\ket\la.\]
Having done this for each $\mu$, we extend semi-linearly to obtain the bar involution on the whole of $\calf^\ba$.  By the above remarks concerning the $\UU$-actions on $\calf^{\ba}$ and $\calf^{\ta}$, this bar involution is compatible with the action of $\UU$ on $\calf^{\ba}$.  In particular, it agrees with the bar involution already defined on $M^\ba$.  We echo the remark of Brundan and Kleshchev \cite[Remark 3.27]{bk} that it would be very interesting to find a construction of this bar involution on $\calf^\ba$ without using twisted Fock spaces.

Once we have defined the bar involution, we can define the canonical basis $\{G^\ba(\mu)\mid \mu\in\calp^r\}$ for $\calf^\ba$.  In fact, the canonical basis element $G^\ba(\mu)$ will be the same as the canonical basis element $G^{\tilde\ba}(\mu)$ for any multicharge $\ta$ with each $\tilde s_k-\tilde s_{k+1}$ large.

We shall need the following dominance property of the canonical basis elements.

\begin{propn}\label{decdom}
Suppose $\mu\in\calp^r$, and write
\[G^\ba(\mu) = \sum_{\la\in\calp^r}d^\ba_{\la\mu}\ket\la.\]
Then $d^\ba_{\la\mu}=0$ unless $\mu\dom\la$.
\end{propn}

\begin{pf}
This follows from \cite[Theorem 2.8 \& Proposition 5.12]{yv1}.
\end{pf}

\sect{A key property of the bar involution}\label{keysec}

In this section we give the details of the construction of the bar involution on a twisted Fock space $\calf^{\ta}$, and prove an important property of the coefficients $b^{\ba}_{\la\mu}$.  Recall that for $\la\in\calp^r$ we define $\tru\la = (\la_2,\dots,\la_r)$; we also define $\tru\ba=(s_2,\dots,s_r)$ for $\ba\in I^r$.

Our aim is to prove the following statement.

\begin{propn}\label{trunc}
Suppose $\ba\in I^r$ for $r>1$ and $\la,\mu\in\calp^r$ with $\mu^{(1)}=\varnothing$.  Then
\[b^{\ba}_{\la\mu} = \begin{cases}
b^{\tru\ba}_{\tru\la\tru\mu} & (\text{if }\la^{(1)}=\varnothing)\\
0 & (\text{otherwise}).
\end{cases}\]
\end{propn}

This gives the following corollary for canonical basis coefficients.

\begin{cory}\label{fcomp}
Suppose $\ba\in I^r$ for $r>1$ and $\mu\in\calp^r$ with $\mu^{(1)}=\varnothing$.  If we write
\begin{align*}
G^{\tru\ba}(\tru\mu) &= \sum_{\nu\in\calp^{r-1}}d^{\tru\ba}_{\nu\tru\mu}\ket\nu,\\
\intertext{then}
G^{\ba}(\mu) &= \sum_{\nu\in\calp^{r-1}}d^{\tru\ba}_{\nu\tru\mu}\ket{\ste\nu}.
\end{align*}
\end{cory}

\begin{pf}
It is straightforward to verify that the vector on the right-hand side of the second equation is bar-invariant, using the bar-invariance of $G^{\tru\ba}(\tru\mu)$ and Proposition \ref{trunc}.  Also, the coefficient of $\ket\mu$ is $d^{\tru\ba}_{\tru\mu\tru\mu}=1$, while all the other coefficients are divisible by $q$.  So by uniqueness of canonical basis elements, this vector must be $G^{\ba}(\mu)$.
\end{pf}

In order to prove Proposition \ref{trunc}, we just need to prove that it holds with $\ba$ replaced by a multicharge $\ta$ for $\ba$ which has $\tilde s_k-\tilde s_{k+1}\gg0$ for each $k$.  To do this, we need to describe in detail how the bar involution on $\calf^{\ta}$ is computed.

Let us define a \emph{wedge} of length $l$ to be a symbol of the form
\[\uu{t_1}\wed\dots\wed\uu{t_l},\]
where $t_1,\dots,t_l\in\bbz$.  We also define a \emph{semi-infinite wedge of charge $s$} to be a symbol
\[\uu{t_1}\wed\uu{t_2}\wed\dots,\]
where $t_1,t_2,\dots\in\bbz$ are such that $t_i=s+1-i$ for $i\gg0$.  We say that a wedge (finite or semi-infinite) is \emph{ordered} if the integers appearing are strictly decreasing.

For fixed $e,r\gs1$, we impose relations on wedges which allow us to express any wedge as a $\bbq(q)$-linear combination of ordered wedges, as follows.  For any integer $t$, let $a(t)\in\inv1e$, $b(t)\in\inv 1r$ and $m(t)\in\bbz$ be such that $t=a(t)+e(b(t)-1)-erm(t)$.  Now given any $t\ls u$ we define $\alpha,\beta$ to be the residues of $a(u)-a(t)$ and $e(b(u)-b(t))$ respectively, modulo $er$.  Then we impose the following relations on wedges of length $2$.
\begin{description}
\item[if $\alpha=\beta=0$:]
\[\uu t\wed\uu u = -\uu u\wed\uu t.\]
\item[if $\alpha>\beta=0$:]
\begin{align*}
\uu t\wed\uu u = &-q^{-1}\uu u\wed\uu t\\
&+(q^{-2}-1)\left(\sum_{m\gs0}q^{-2m}\uu{u-\alpha-erm}\wed\uu{t+\alpha+erm}-\sum_{m\gs1}q^{1-2m}\uu{u-erm}\wed\uu{t+erm}\right).
\end{align*}
\item[if $\beta>\alpha=0$:]
\begin{align*}
\uu t\wed\uu u = &-q\uu u\wed\uu t\\
&+(q^2-1)\left(\sum_{m\gs0}q^{2m}\uu{u-\beta-erm}\wed\uu{t+\beta+erm}-\sum_{m\gs1}q^{2m-1}\uu{u-erm}\wed\uu{t+erm}\right).
\end{align*}
\item[if $\alpha,\beta>0$:]
\begin{align*}
\uu t\wed\uu u = &-\uu u\wed\uu t\\
&+(q-q^{-1})\sum_{m\gs0}\frac{q^{2m+1}+q^{-2m-1}}{q+q^{-1}}\uu{u-\beta-erm}\wed\uu{t+\beta+erm}\\
&+(q^{-1}-q)\sum_{m\gs0}\frac{q^{2m+1}+q^{-2m-1}}{q+q^{-1}}\uu{u-\alpha-erm}\wed\uu{t+\alpha+erm}\\
&+(q-q^{-1})\sum_{m\gs0}\frac{q^{2m+2}-q^{-2m-2}}{q+q^{-1}}\uu{u-\alpha-\beta-erm}\wed\uu{t+\alpha+\beta+erm}\\
&+(q^{-1}-q)\sum_{m\gs1}\frac{q^{2m}-q^{-2m}}{q+q^{-1}}\uu{u-erm}\wed\uu{t+erm}.
\end{align*}
\end{description}
In each of these expressions, the summation continues only as long as the wedges are ordered.

For $l>2$, or for semi-infinite wedges, the ordering relations are defined by imposing the above relations in each adjacent pair of positions.  Now we define the $l$-wedge space to be the $\bbq(q)$-vector space spanned by all wedges of length $l$, modulo the ordering relations.  We also define the semi-infinite wedge space of charge $s$ to be the $\bbq(q)$-vector space spanned by the set of all semi-infinite wedges of charge $s$, modulo the ordering relations.  In each of these spaces, the set of ordered wedges is a basis.

In order to avoid ambiguity when comparing different values of $r$, we may decorate the wedge symbol $\wed$ as $\swed r$ to indicate the particular value of $r$ used in the straightening relations.

The construction of the bar involution relies on encoding a pair $(\la,\ta)$ (where $\la\in\calp^r$ and $\ta$ is a multicharge) as an ordered wedge.  We set $s=\tilde s_1+\dots+\tilde s_r$, and define a semi-infinite wedge of charge $s$ as follows.\label{wedgedef}

For each $k\in\inv 1r$ and each $i\gs1$ set
\[\beta^{(k)}_i = \la^{(k)}_i+\tilde s_k+1-i.\]
Write this integer in the form
\[\beta^{(k)}_i = a-em\]
with $m\in\bbz$ and $a\in\inv1e$, and then set
\[\check\beta^{(k)}_i = a+e(k-1)-erm.\]
The integers $\check\beta^{(k)}_i$ for $k\in\inv1r$ and $i\gs1$ are distinct and bounded above, so we may arrange them in strictly decreasing order as $t_1>t_2>\dots$.  Then the ordered wedge $\kt{\la,\ta}$ corresponding to $\la$ and $\ta$ is
\[\uu{t_1}\wed\uu{t_2}\wed\dots.\]
It is easy to check that this is a semi-infinite wedge of charge $s$.  Conversely, each ordered semi-infinite wedge of charge $s$ is equal to $\kt{\la,\ta}$ for some $r$-partition $\la$ and some multicharge $\ta$ with sum $s$.

Now we can define the bar involution on $\calf^{\ta}$.  Given an $r$-multipartition $\mu$, we write
\[\kt{\mu,\ta} = \uu{t_1}\wed\uu{t_2}\wed\dots;\]
we choose $l\gg0$, and set
\[\tk{\mu,\ta} = \uu{t_l}\wed\uu{t_{l-1}}\wed\dots\wed\uu{t_1}\wed\uu{t_{l+1}}\wed\uu{t_{l+2}}\wed\dots.\]
Using the ordering relations, we express $\tk{\mu,\ta}$ as a linear combination of ordered wedges.  It is easy to show (by considering residues modulo $er$) that each of the ordered wedges that occurs has the form $\kt{\la,\ta}$ for some $r$-multipartition $\la$, i.e.\ we have a finite sum
\[\tk{\mu,\ta} = \sum_\la c_{\la\mu}\kt{\la,\ta}\]
with each $c_{\la\mu}\in\bbq(q)$.  Moreover, the coefficient $c_{\mu\mu}$ is non-zero, so we can define
\[\ol{\ket\mu} = \sum_\la \frac{c_{\la\mu}}{c_{\mu\mu}}\ket\la.\]
This definition is independent of $l$, provided we choose $l$ sufficiently large.  This defines the bar involution on the basis elements $\ket\mu$, and we extend semi-linearly to obtain the bar involution on the whole of $\calf^{\ta}$.

Now we set about proving Proposition \ref{trunc}.  The calculations used here are similar to those used in \cite{runner}, though actually rather simpler.  Since there is nothing to prove when $r=1$, we assume for the rest of this section that $r\gs2$.

Recalling the definition of $b(t)$ for $t\in\bbz$ from above, we define $\one = b^{-1}(1)$.  In other words, $\one$ consists of all integers whose residue modulo $er$ lies in $\inv1e$.  Now we define a map $\psi:\bbz\setminus\one\to\bbz$: given $t\in\bbz\setminus\one$, we define $a(t),b(t),m(t)$ as above, and set
\[\psi(t) = a(t)+e(b(t)-2)-e(r-1)m(t).\]
Then $\psi$ is an order-preserving bijection from $\bbz\setminus\one$ to $\bbz$.  Furthermore, the following relationship is easy to check from the straightening relations.

\begin{lemma}\label{rrm1}
Suppose $t_{ij}\in\bbz\setminus\one$ and $c_i\in\bbq(q)$ for $1\ls i\ls m$ and $1\ls j\ls l$.  Then
\[\sum_{i=1}^mc_i\uu{t_{i1}}\swed r\dots\swed r\uu{t_{il}}=0\]
if and only if
\[\sum_{i=1}^mc_i\uu{\psi(t_{i1})}\swed {r-1}\dots\swed {r-1}\uu{\psi(t_{il})}=0.\]
\end{lemma}

Using $\psi$, we can describe the relationship between the wedges $\kt{\la,\ta}$ and $\kt{\tru\la,\tru\ta}$.

\begin{lemma}\label{compareweds}
Suppose $\mu\in\calp^r$ and $\ta$ is a multicharge, and write
\[\kt{\mu,\ta} = \uu{t_1}\wed\uu{t_2}\wed\dots.\]
Then, if we write the elements of $\{t_1,t_2,\dots\}\setminus\one$ as $u_1>u_2>\dots$, we have
\[\kt{\tru\mu,\tru\ta} = \uu{\psi(u_1)}\wed\uu{\psi(u_2)}\wed\dots.\]
Furthermore, if $\mu^{(1)}=\varnothing$ and $\tilde s_1\gg\tilde s_k$ for each $k\in\inv2r$, then there is an integer $d\gs t_1$ such that
\[\{t_1,t_2,\dots\}\cap\one = \bbz_{\ls d}\cap\one.\]
\end{lemma}

\begin{pf}
This is easy to check from the definition of $\kt{\mu,\ta}$.
\end{pf}

Lemma \ref{compareweds} allows us to compare the computations of $\ol{s_\mu}$ and $\ol{s_{\tru\mu}}$ (in the twisted Fock spaces $\calf^{\ta}$ and $\calf^{\tru\ta}$, respectively) when $\mu^{(1)}=\varnothing$ and $\tilde s_1\gg\dots\gg\tilde s_r$.  The idea is that we write
\[\kt{\mu,\ta} = \uu{t_1}\wed\uu{t_2}\wed\dots,\]
and then straighten the finite wedge
\[\uu{t_l}\wed\dots\wed\uu{t_1}\]
for suitably large $l$.  We do this by first moving the terms $\uu{t_j}$ with $t_j\in\one$ to the beginning; then we order these terms, and we separately order the remaining terms (employing Lemma \ref{rrm1}).  Finally, we recombine terms to obtain a linear combination of ordered wedges.  In the next few results, we check the details of this procedure.

\begin{lemma}\label{portm}
Suppose $c\ls d$ are integers and $t_1,\dots,t_l\in\inv cd$, and let
\[W=\uu{t_1}\wed\dots\wed\uu{t_l}.\]
When we express $W$ as a linear combination or ordered wedges using the straightening relations, each ordered wedge $\uu{u_1}\wed\dots\wed\uu{u_l}$ 
that occurs satisfies
\[u_1,\dots,u_l\in\inv cd\]
and
\[\left|\left\{\left.j\in\inv1l\ \right|\ u_j\in\one\right\}\right|=\left|\left\{\left.j\in\inv1l\ \right|\ t_j\in\one\right\}\right|.\]
\end{lemma}

\begin{pf}
This is easy to check from the straightening relations.
\end{pf}

\begin{cory}\label{aminus}
Suppose $t_1,\dots,t_l\in\inv cd$.  Suppose that
\[\left|\left\{\left.j\in\inv1l\ \right|\ t_j\in\one\right\}\right|>\big|\inv cd\cap\one\big|.\]
Then the wedge $W=\uu{t_1}\wed\dots\wed\uu{t_l}$ equals zero.
\end{cory}

\begin{pf}
$W$ can be written as a linear combination of ordered wedges using the straightening relations, and each ordered wedge $\uu{u_1}\wed\dots\wed\uu{u_l}$ occurring must satisfy the conditions in Lemma \ref{portm}.  But the hypotheses on $t_1,\dots,t_l$ mean that there are no such ordered wedges, and so $w$ must equal zero.
\end{pf}

\begin{lemma}\label{a}
Suppose $c\ls d$ are integers, and write the elements of the set $\inv cd\cap\one$ as $u_1>\dots>u_n$.  Suppose $v_1>\dots>v_m$ are elements of $\inv cd\setminus\one$, and label the elements of the set $\{v_1,\dots,v_m,u_1,\dots,u_n\}$ in decreasing order as $t_1>\dots>t_{m+n}$.  Then the wedge
\[W = \uu{t_{m+n}}\wed\uu{t_{m+n-1}}\wed\dots\wed\uu{t_1}\]
is equal to a scalar multiple of the wedge
\[W' = \uu{u_n}\wed\dots\wed\uu{u_1}\wed\uu{v_m}\wed\dots\wed\uu{v_1}.\]
\end{lemma}

\begin{pf}
If $t_{m+n}=u_n$, then $W$ and $W'$ have the same first term; by induction on $n$ (replacing $c$ with $t_{m+n}+1$) the wedge obtained by removing the first term from $W$ is proportional to the wedge obtained by removing the first term from $W'$, so $W$ and $W'$ are proportional too.  So we may assume that $t_{m+n}=v_m$.  We also assume that $m=1$; the general case follows by induction on $m$.  So by assumption we have
\[W=\uu{v_1}\wed\uu{u_n}\wed\dots\wed\uu{u_1},\qquad W'=\uu{u_n}\wed\dots\wed\uu{u_1}\wed\uu{v_1}.\]

Using induction on $n$ again (replacing $d$ with $u_1-1$) $W$ is equal to a multiple of
\[W_1=\uu{u_n}\wed\dots\wed\uu{u_2}\wed\uu{v_1}\wed\uu{u_1};\]
applying the straightening relations to $\uu{v_1}\wed\uu{u_1}$, we find that $W_1$ equals a scalar multiple of $W'$ plus a linear combination of wedges of the form
\[\uu{u_n}\wed\dots\wed\uu{u_2}\wed\uu{w}\wed\uu{w'}\]
in which $w,w'$ lie strictly between $v_1$ and $u_1$, and one of $w,w'$ lies in $\one$.  Now by Corollary \ref{aminus} (with $u_1-1$ in place of $d$) each such wedge is equal to zero, so $W_1$ is proportional to $W'$.
\end{pf}

\begin{lemma}\label{b}
Suppose $c\ls d$ are integers, and write the elements of the set $\inv cd\cap\one$ as $u_1>\dots>u_n$.  Then the wedge
\[W=\uu{u_n}\wed\dots\wed\uu{u_1}\]
is equal to a scalar multiple of
\[W'=\uu{u_1}\wed\dots\wed\uu{u_n}.\]
\end{lemma}

\begin{pf}
When we write $W$ as a linear combination of ordered wedges using the straightening relations, each ordered wedge that occurs satisfies the conditions in Lemma \ref{portm}.  But the only such ordered wedge is $W'$.
\end{pf}

Now given integers $c\ls v$, define
\[X_c(v) = \big|\inv cv\cap\one\big|,\qquad Y_c(v) = \big|\inv cv\cap\one\cap(v+e\bbz)\big|.\]

\begin{lemma}\label{samexy}
Suppose $c\ls d$ and $v_1,\dots,v_m\in\inv cd\setminus\one$.  When the wedge
\[\uu{v_1}\wed\dots\wed\uu{v_m}\]
is written as a linear combination of ordered wedges using the straightening relations, each wedge $\uu{w_1}\wed\dots\wed\uu{w_m}$ that occurs with non-zero coefficient satisfies
\[\sum_{i=1}^mX_c(w_i)=\sum_{i=1}^mX_c(v_i),\qquad
\sum_{i=1}^mY_c(w_i)=\sum_{i=1}^mY_c(v_i).\]
\end{lemma}

\begin{pf}
Consider the case $m=2$.  In this case, $\uu{v_1}\wed\uu{v_2}$ is equal to a linear combination of ordered wedges $\uu{w_1}\wed\uu{w_2}$ for which (recalling the functions $a,b$ from above):
\begin{itemize}
\item
$w_1,w_2\in\inv cd$;
\item
$w_1+w_2=v_1+v_2$;
\item
$\left\{a(w_1),a(w_2)\right\}=\left\{a(v_1),a(v_2)\right\}$;
\item
$\left\{b(w_1),b(w_2)\right\}=\left\{b(v_1),b(v_2)\right\}$.
\end{itemize}
From these properties, it follows easily that $X_c(w_1)+X_c(w_2)=X_c(v_1)+X_c(v_2)$ and $Y_c(w_1)+Y_c(w_2)=Y_c(v_1)+Y_c(v_2)$.

The case $m>2$ follows by applying the above case each time a straightening rule is applied.
\end{pf}

\begin{lemma}\label{releave}
Suppose $c\ls d$, and write the elements of $\inv cd\cap\one$ as $u_1>\dots>u_n$.  Suppose $v_1>\dots>v_m$ are elements of $\inv cd\setminus\one$, and label the elements of the set $\{v_1,\dots,v_m,u_1,\dots,u_n\}$ as $t_1>\dots>t_{m+n}$.  Then the wedge
\[W=\uu{u_1}\wed\dots\wed\uu{u_n}\wed\uu{v_1}\wed\dots\wed\uu{v_m}\]
equals
\[\left(\prod_{i=1}^m(-1)^{X_c(v_i)}q^{Y_c(v_i)}\right)\uu{t_1}\wed\dots\wed\uu{t_{m+n}}.\]
\end{lemma}

\begin{pf}
We may assume that $t_1=v_1$; otherwise, the result follows by induction on $n$ (replacing $d$ with $u_1-1$).  We can also assume (using induction on $m$) that $m=1$.  So we assume that
\[W=\uu{u_1}\wed\dots\wed\uu{u_n}\wed\uu{v_1}\]
with $v_1>u_1$, and we want to show that
\[W = (-1)^nq^Y\uu{v_1}\wed\uu{u_1}\wed\dots\wed\uu{u_n},\]
where $Y=\left|\{u_1,\dots,u_n\}\cap (v_1+e\bbz)\right|$.

Applying the straightening rules to $\uu{u_n}\wed\uu{v_1}$, we find that
\[W=-q'\uu{u_1}\wed\dots\wed\uu{u_{n-1}}\wed\uu{v_1}\wed\uu{u_n}\]
where $q'$ equals $q$ if $u_n$ and $v_1$ are congruent modulo $e$, and $1$ otherwise.  (The other terms arising from applying the straightening relation vanish by Corollary \ref{aminus}.)

Now induction on $n$ (replacing $c$ with $u_n+1$) gives the result.
\end{pf}

\begin{pf}{Proof of Proposition \ref{trunc}}
It suffices to prove the result with $\ba$ replaced by a multicharge $\ta$ for $\ba$ such that $\tilde s_k-\tilde s_{k+1}\gg0$ for each $k$.  So we choose such a multicharge, and write
\[\kt{\mu,\ta} = \uu{t_1}\wed\uu{t_2}\wed\dots.\]
We write the elements of $\{t_1,t_2,\dots\}\cap\one$ as $u_1>u_2>\dots$, and the elements of $\{t_1,t_2,\dots\}\setminus\one$ as $v_1>v_2>\dots$.  By Lemma \ref{compareweds}, the set $\{u_1,u_2,\dots\}$ consists of all elements of $\one$ which are less than or equal to $u_1$.

To compute the effect of the bar involution on $s_\mu$, we straighten the wedge
\[W = \uu{t_l}\wed\dots\wed\uu{t_1},\]
where $l\gg0$ is fixed.  Let $m,n$ be such that
\[\{t_1,\dots,t_l\}=\{v_1,\dots,v_m\}\cup\{u_1,\dots,u_n\};\]
if we put $c=t_l$, $d=t_1$, then we have $\{u_1,\dots,u_n\} = \inv cd\cap\one$ and $v_1,\dots,v_m\in\inv cd\setminus\one$, so by Lemma \ref{a} and Lemma \ref{b} $W$ is equal to a scalar multiple of
\[\uu{u_1}\wed\dots\wed\uu{u_n}\wed\uu{v_m}\wed\dots\wed\uu{v_1}.\]
Now write the wedge $\uu{v_m}\wed\dots\wed\uu{v_1}$ as a linear combination of ordered wedges:
\[\uu{v_m}\wed\dots\wed\uu{v_1} = \sum_{i=1}^N\alpha_i\uu{v^i_1}\wed\dots\wed\uu{v^i_m}.\tag*{($\ast$)}\]
For each $i$, let $t^i_1,\dots,t^i_l$ be the sequence obtained by putting the integers $v^i_1,\dots,v^i_m,u_1,\dots,u_n$ in decreasing order.  Then by Lemma \ref{samexy} and Lemma \ref{releave} we find that $W$ is equal to a scalar multiple of
\[\sum_{i=1}^N\alpha_i\uu{t^i_1}\wed\dots\wed\uu{t^i_l}.\tag*{($\dagger$)}\]

Now we consider how to compute $\ol{\ket{\tru\mu}}$ in $\calf^{\tru\ta}$.  From Lemma \ref{compareweds}, we have
\[\kt{\tru\mu,\tru\ta} = \uu{\psi(v_1)}\wed\uu{\psi(v_2)}\wed\dots.\]
By Lemma \ref{rrm1} and ($\ast$) we have
\[\uu{\psi(v_m)}\swed {r-1}\dots\swed {r-1}\uu{\psi(v_1)} = \sum_{i=1}^N\alpha_i\uu{\psi(v^i_1)}\swed {r-1}\dots\swed {r-1}\uu{\psi(v^i_m)}.\]
Since $l$ (and hence $m$) is large, we therefore find that for each $i$ there is a multipartition $\nu(i)\in\calp^{r-1}$ such that
\[\kt{\nu(i),\tru\ta}=\uu{\psi(v^i_1)}\wed\dots\wed\uu{\psi(v^i_m)}\wed\uu{\psi(v_{m+1})}\wed\uu{\psi(v_{m+2})}\wed\dots.\]
Hence there is $\alpha\in\bbq(q)$ (independent of $i$) such that $\alpha_i = \alpha b^{\tru\ta}_{\nu(i)\mu}$ for each $i$; moreover, each $\nu$ for which $b^{\tru\ta}_{\nu\mu}\neq0$ occurs as some $\nu(i)$.

Now by Lemma \ref{compareweds} we have
\[\uu{t^i_1}\wed\dots\wed\uu{t^i_l}\wed\uu{t_{l+1}}\wed\uu{t_{l+2}}\wed\dots = \kt{\nu(i)_+,\ta},\]
and the result follows from ($\dagger$).
\end{pf}

\sect{An LLT-type algorithm}\label{algsec}

Now we can give our algorithm which generalises the LLT algorithm.  As mentioned above, our algorithm actually computes the canonical basis of $\MM{\ba}\cong M^{(s_1)}\otimes\dots\otimes M^{(s_r)}$.

Since the canonical basis elements $G^{(s_k)}(\mu)$ indexed by $e$-regular partitions $\mu$ form a basis for $M^{(s_k)}$, the tensor product $M^{(s_1)}\otimes\dots\otimes M^{(s_r)}$ has a basis consisting of all vectors $G^{(s_1)}(\mu^{(1)})\otimes\dots\otimes G^{(s_r)}(\mu^{(r)})$, where $\mu^{(1)},\dots,\mu^{(r)}$ are $e$-regular partitions.  Translating this to the Fock space $\calf^\ba$, we find that $\MM{\ba}$ has a basis consisting of vectors
\[\GG^\ba(\mu) = \sum_{\la\in\calp^r}d^{(s_1)}_{\la^{(1)}\mu^{(1)}}\dots d^{(s_r)}_{\la^{(r)}\mu^{(r)}}\ket\la\]
for all $e$-multiregular multipartitions $\mu$.  In fact, we want to show that the canonical basis vectors $G^\ba(\mu)$ for $e$-multiregular $\mu$ form a basis for $\MM{\ba}$; this implies in particular that the span of these vectors is a $U$-submodule of $\calf^\ba$, which will enable our recursive algorithm to work.

\begin{lemma}\label{multiregutil}
Suppose $G = \sum_{\la\in\calp^r}g_\la\ket\la\in\MM{\ba}$.  If $\la\in\calp^r$ is such that $g_\la\notin q\bbq[q]$, then there is an $e$-multiregular multipartition $\nu$ such that $\nu\dom\la$, $|\nu^{(k)}|=|\la^{(k)}|$ for all $k$ and $g_\nu\notin q\bbq[q]$.
\end{lemma}

\begin{pf}
We may write $G$ as
\[G = \sum_{\nu\in\reg^r}h_\nu\GG^\ba(\nu)\]
with $h_\nu\in\bbq(q)$ for each $\nu$; then we have
\[g_\la = \sum_{\nu\in\reg^r} h_\nu d^{(s_1)}_{\la^{(1)}\nu^{(1)}}\dots d^{(s_r)}_{\la^{(r)}\nu^{(r)}}.\]
Since $d^{(s_k)}_{\la^{(k)}\nu^{(k)}}$ can be non-zero only if $|\la^{(k)}|=|\nu^{(k)}|$ and $\nu^{(k)}\dom\la^{(k)}$, we may restrict the range of summation to only those $\nu$ which have $|\la^{(k)}|=|\nu^{(k)}|$ and $\nu^{(k)}\dom\la^{(k)}$ for all $k$ (and hence $\nu\dom\la$).

Since $g_\la\notin q\bbq[q]$ but each $d^{(s_k)}_{\la^{(k)}\nu^{(k)}}$ is a polynomial in $q$, we must have $h_\nu\notin q\bbq[q]$ for some $\nu$.  If we choose such a $\nu$ which is maximal with respect to the dominance ordering, then we have
\[g_\nu-h_\nu=\sum_{\nu\domsby\xi\in\reg^r}h_\xi d^{(s_1)}_{\nu^{(1)}\xi^{(1)}}\dots d^{(s_r)}_{\nu^{(r)}\xi^{(r)}}\in q\bbq[q].\]
Now the fact that $h_\nu\notin q\bbq[q]$ implies that $g_\nu\notin q\bbq[q]$.
\end{pf}

Now we can deduce the following.

\begin{propn}\label{algo}
The canonical basis vectors $G^\ba(\mu)$ indexed by $e$-multiregular $\mu$ form a basis for the module $\MM{\ba}$.
\end{propn}

\begin{pf}
All we need to do is show that $G^\ba(\mu)$ lies in $\MM{\ba}$ for each $e$-multiregular $\mu$; since the canonical basis vectors are linearly independent and $G^\ba(\mu)$ and $\GG^\ba(\mu)$ are both weight vectors of the same weight, the result follows by considering the dimensions of weight spaces in $\MM{\ba}$.

We prove that $G^\ba(\mu)$ lies in $\MM{\ba}$ for each $\mu\in\reg^r$ by induction on $r$ and, for fixed $r$, by induction on $|\mu^{(1)}|$.  When $r=1$ there is nothing to prove, since then $G^\ba(\mu)=\GG^\ba(\mu)$.

Suppose $r>1$ and $\mu^{(1)}=\varnothing$.  By induction on $r$, $G^{\tru\ba}(\tru\mu)$ can be written as a linear combination
\begin{align*}
G^{\tru\ba}(\tru\mu) &= \sum_{\nu\in\reg^{r-1}}c_\nu\GG^{\tru\ba}(\nu).\\
\intertext{By Corollary \ref{fcomp} and the fact that $G^{(s_1)}(\varnothing)=s_{\varnothing}$, we therefore have}
G^\ba(\mu) &= \sum_{\nu\in\reg^{r-1}}c_\nu\GG^\ba(\ste\nu),
\end{align*}
so $G^\ba(\mu)\in\MM{\ba}$.

Now consider the case where $r>1$ and $|\mu^{(1)}|>0$.  Using the LLT algorithm, we can write $G^{(s_1)}(\mu^{(1)})$ as $us_\varnothing$ in the Fock space $\calf^{(s_1)}$, for some $u\in\UU$.  Moreover, we can choose $u$ to be a $\bbz[q+q^{-1}]$-linear combination of products of divided powers $f_i^{(a)}$.

Recall that we write $\mu_0$ to mean $(\varnothing,\mu^{(2)},\dots,\mu^{(r)})$.  By induction on $|\mu^{(1)}|$, we have $G^\ba(\mu_0)\in\MM{\ba}$.  So if we define $G=uG^\ba(\mu_0)$, then since $\MM{\ba}$ is a $\UU$-submodule of $\calf^\ba$, we have $G\in\MM{\ba}$.  Furthermore, because $G^\ba(\mu_0)$ is bar-invariant and because of the properties of the element $u$, $G$ is also bar-invariant.  If we write $G=\sum_{\la\in\calp^r}g_\la\ket\la$, then, using the rule for the actions of the $f_i$, we find that
\begin{alignat*}2
g_\la &= d^{(s_1)}_{\la^{(1)}\mu^{(1)}}d^\ba_{\la_0\mu_0}\qquad&&\text{if }|\la^{(1)}|=|\mu^{(1)}|,\\
\intertext{while}
g_\la &= 0&&\text{if }|\la^{(1)}|>|\mu^{(1)}|.
\end{alignat*}
In particular, if $|\la^{(1)}|\gs|\mu^{(1)}|$ and $\la\neq\mu$, then $q\mid g_\la$ ($\ddagger$).  Furthermore, each $g_\la$ lies in $\bbz[q,q^{-1}]$.

Now consider the expansion of $G$ as a linear combination of canonical basis elements.  Since $g_\mu=1$ and $g_\nu=0$ for any $\nu\doms\mu$, Proposition \ref{decdom} implies that $G$ equals $G^\ba(\mu)$ plus a linear combination of canonical basis elements $G^\ba(\nu)$ with $\nu\ndom\mu$; because $G^\ba(\mu)$ is bar-invariant, the coefficients of these canonical basis elements all lie in $\bbz[q+q^{-1}]$.  This means that one can apply the same procedure as in the LLT algorithm to `strip off' the terms $G^\ba(\nu)$ with $\nu\neq\mu$ and recover $G^\ba(\mu)$.  This is done as follows:
\begin{itemize}
\item
if there is no $\nu\neq\mu$ such that $g_\nu\notin q\bbz[q]$, then stop.
\item
otherwise, choose such a $\nu$ which is maximal with respect to the dominance ordering, and let $\alpha$ be the unique element of $\bbz[q+q^{-1}]$ such that $g_\nu-\alpha\in q\bbz[q]$.  Replace $G$ with $G-\alpha G^\ba(\nu)$, and repeat.
\end{itemize}
At each stage, the vector $G$ lies in $\MM{\ba}$, and so by Lemma \ref{multiregutil} and ($\ddagger$), the multipartition $\nu$ involved must be $e$-multiregular and must satisfy $|\nu^{(1)}|<|\mu^{(1)}|$.  Therefore by induction we have $G^\ba(\nu)\in\MM{\ba}$, and so the new vector $G-\alpha G^\ba(\nu)$ lies in $\MM{\ba}$.

At the end of this procedure, we are left with the canonical basis vector $G^\ba(\mu)$, and this lies in $\MM{\ba}$.
\end{pf}

Proposition \ref{algo} enables us to construct canonical basis vectors labelled by $e$-multiregular multipartitions recursively.  As in the LLT algorithm, the idea is that to construct the canonical basis vector $G^\ba(\mu)$, we construct an auxiliary vector $A(\mu)$ which is bar-invariant, and which we know equals $G^\ba(\mu)$ plus a linear combination of `lower' canonical basis vectors; the bar-invariance of $A(\mu)$, together with dominance properties, allows these lower terms to be stripped off.  In our algorithm, we take additional care to make sure that $A(\mu)$ lies in $\MM{\ba}$; then we know by Proposition \ref{algo} that all the canonical basis vectors occurring in $A(\mu)$ are labelled by $e$-multiregular multipartitions, and therefore we can assume that these have already been constructed.

In fact, the proof of Proposition \ref{algo}, combined with the construction in the LLT algorithm, gives us our algorithm.  We formalise this as follows.

Our algorithm is recursive, using a partial order on multipartitions which is finer than the dominance order: define $\mu\cgs\nu$ if either $|\mu^{(1)}|>|\nu^{(1)}|$ or $\mu^{(1)}\dom\nu^{(1)}$.  We assume when computing $G^\ba(\mu)$ for $\mu\in\reg^r$ that we have already computed the vector $G^{\tru\ba}(\tru\mu)$, and that we have computed $G^\ba(\nu)$ for all $\nu\in\reg^r$ with $\mu\cg\nu$.

\begin{enumerate}
\item
If $\mu=\varnothing^r$, then $G^\ba(\mu)=\ket{\varnothing^r}$.
\item
If $\mu\neq\varnothing^r$ but $\mu^{(1)}=\varnothing$, then compute the canonical basis vector $G^{\tru\ba}(\tru\mu)$.  Then $G^\ba(\mu)$ is given by
\[G^\ba(\mu) = \sum_{\nu\in\calp^{r-1}}d^{\tru\ba}_{\nu\tru\mu}\ket{\ste\nu}.\]
\item
If $\mu^{(1)}\neq\varnothing$, then apply the following procedure.
\begin{enumerate}
\item
Let $\mu_0=(\varnothing,\mu^{(2)},\dots,\mu^{(r)})$, and compute $G^\ba(\mu_0)$.
\item
Let $a_1,\dots,a_t$ be the sizes of the non-empty ladders of $\mu^{(1)}$, and $i_1,\dots,i_t$ their residues.  Define $A = f_{i_t}^{(a_t)}\dots f_{i_1}^{(a_1)}G^\ba(\mu_0)$. Write $A=\sum_{\nu\in\calp^r}a_\nu\ket\nu$.  
\item\label{stepc}
If there is no $\nu\neq\mu$ for which $a_\nu\notin q\bbz[q]$, then stop.  Otherwise, take such a $\nu$ which is maximal with respect to the dominance order, let $\alpha$ be the unique element of $\bbz[q+q^{-1}]$ for which $a_\nu-\alpha\in q\bbz[q]$, replace $A$ by $A-\alpha G^\ba(\nu)$, and repeat.  The remaining vector will be $G^\ba(\mu)$.
\end{enumerate}
\end{enumerate}

The vector $A$ computed in step 3 is a bar-invariant element of $\MM{\ba}$, because $G^\ba(\mu_0)$ is.  Hence by Proposition \ref{algo} $A$ is a $\bbq(q+q^{-1})$-linear combination of canonical basis vectors $G^\ba(\nu)$ with $\nu\in\reg^r$.  Furthermore, the rule for applying $f_i$ to a multipartition and the combinatorial results used in the LLT algorithm imply that $a_\mu=1$, and that if $a_\la\neq0$, then $\mu\cgs\la$.  
In particular, the partition $\nu$ appearing in step 3(c) satisfies $\mu\cg\nu$; moreover, when $\alpha G^\ba(\nu)$ is subtracted from $A$, the condition that $a_\mu=1$ and $a_{\la}$ is non-zero only for $\mu\cgs\la$ remains true (because of Proposition \ref{decdom} and the fact that the dominance order refines the order $\cgs$).  So we can repeat, and complete step 3(c).


\sect{An example and further remarks}\label{egsec}

\subs{An example}

Let us take $e=r=2$, and write the set $I=\bbz/2\bbz$ as $\{0,1\}$.  Take $\ba=(0,0)$.
\begin{itemize}
\item
First let us compute the canonical basis element $G^\ba(((2,1),(1)))$.  In the level $1$ Fock space $\calf^{(0)}$, we have $G^{(0)}((1))=\ket{(1)}$ (where the partition $(1)$ really stands for the $1$-multipartition $((1))$).  The non-empty ladders of the partition $(2,1)$ are $\lad_1$ and $\lad_2$, of lengths $1,2$ and residues $0,1$ respectively.  So we compute
\[A=f_1^{(2)}f_0\ket{(\varnothing,(1))} = \ket{((2,1),(1))}+q\ket{((2),(2))}+q^2\ket{((2),(1^2))}+q^2\ket{((1^2),(2))}+q^3\ket{((1^2),(1^2))}+q^4\ket{((1),(2,1))}.\]
Since the coefficients in $A$ (apart from the leading one) are divisible by $q$, we have $A=G^\ba(((2,1),(1)))$.
\item
Next we compute $G^\ba(((4),\varnothing))$.  This time our auxiliary vector is
\begin{align*}
A=f_1f_0f_1f_0\ket{(\varnothing,\varnothing)} = &\ \ket{((4),\varnothing)}+q\ket{((3,1),\varnothing)}+q\ket{((2,1^2),\varnothing)}+q^2\ket{((1^4),\varnothing)}+(1+q^2)\ket{((2,1),(1))}\\
&+2q\ket{((2),(2))}+2q^2\ket{((2),(1^2))}+2q^2\ket{((1^2),(2))}+2q^3\ket{((1^2),(1^2))}+(q^2+q^4)\ket{((1),(2,1))}\\
&+q^2\ket{(\varnothing,(4))}+q^3\ket{(\varnothing,(3,1))}+q^3\ket{(\varnothing,(2,1^2))}+q^4\ket{(\varnothing,(1^4))}.
\end{align*}
And so we have
\begin{align*}
G^\ba(((4),\varnothing)) = A-G^\ba(((2,1),(1))) = &\ \ket{((4),\varnothing)}+q\ket{((3,1),\varnothing)}+q\ket{((2,1^2),\varnothing)}+q^2\ket{((1^4),\varnothing)}+q^2\ket{((2,1),(1))}\\
&+q\ket{((2),(2))}+q^2\ket{((2),(1^2))}+q^2\ket{((1^2),(2))}+q^3\ket{((1^2),(1^2))}+q^2\ket{((1),(2,1))}\\
&+q^2\ket{(\varnothing,(4))}+q^3\ket{(\varnothing,(3,1))}+q^3\ket{(\varnothing,(2,1^2))}+q^4\ket{(\varnothing,(1^4))}.
\end{align*}
\end{itemize}

\subs{The case $e=\infty$}\label{infsec}

In this section, we indicate very briefly how our results can be extended to the case $e=\infty$.  Normally in this subject, this extension is straightforward: $U_q(\widehat{\mathfrak{sl}}_e)$ must be replaced with $U_q(\mathfrak{sl}_\infty)$, $\zez$ is replaced with $\bbz$, the set of $e$-regular partitions is replaced with the set of all partitions, and other definitions and results are modified appropriately.  However, for the subject matter of this paper, the situation is more complicated, because the definition of Uglov's twisted Fock spaces does not work in the case $e=\infty$; so a little more discussion is merited.

To get around the difficulty of not having a twisted Fock space, one can `approximate' the case $e=\infty$ using a value of $e$ which is large relative to the partitions in question.  Formally, one restricts attention to multipartitions of size at most $n$, by regarding the Fock space just as a module for the negative part $\UU^-$ of $\UU$ and then passing to the quotient $\calf^\ba_{\ls n}$ by the submodule spanned by all $\ket\la$ with $|\la|>n$.  Now given $\ba\in\bbz^r$, one can take a value of $e$ which is large relative to $n$ and $\ba$, and define the bar involution on $\ket\la$ for $|\la|\ls n$ by using the bar involution on $\calf^{\ba+e\bbz}$ (where $\ba+e\bbz$ means $(s_1+e\bbz,\dots,s_r+e\bbz)$).  Because $e$ is large, the actions of $U_q^-(\widehat{\mathfrak{sl}}_e)$ and $\UU^-$ `agree' on $\calf_{\ls n}$, so this bar involution is compatible with the action of $\UU^-$.

One does this for all $n$, and then defines a bar involution on the whole of $\calf^\ba$ by taking a limit.  Of course, one needs to check that this construction of the bar involution on $\calf_{\ls n}^\ba$ is independent of the choice of $e\gg0$.  This is not too difficult to show using the straightening relations, but in fact we can show this using our algorithm for computing canonical basis elements.  Consider applying our algorithm to compute the canonical basis elements $G^{\ba+e\bbz}(\mu)$ for all $e$-multiregular multipartitions $\mu$ with $|\mu|\ls n$.  The crucial point is that when $e$ is very large,
\begin{itemize}
\item
the implementation of the algorithm doesn't actually involve the integer $e$ in a non-trivial way, and
\item
\emph{every} multipartition of size at most $n$ is $e$-multiregular.
\end{itemize}
Hence one can actually compute the canonical basis for the whole of the truncated Fock space $\calf^{\ba}_{\ls n}$ by this algorithm, and this basis is independent of $e$; since one can recover the bar involution from the canonical basis, this means that the bar involution on $\calf^\ba_{\ls n}$ is independent of the choice of large $e$.

\begin{eg}
Let $r=3$ and $\ba=(0,1,0)$, and take $\mu=((2,1),\varnothing,(1))$.  We shall compute $G^{\ba}(\mu)$ by computing $G^{\hat{\ba}}(\mu)$ for $e\gs4$; for this example write $n+e\bbz$ as $\hat n$, for any $n\in\bbz$.  Starting with the level $1$ Fock space $\calf^{(\hat0)}$, it is easy to compute $G^{(\hat0)}((1))=\ket{(1)}$.  Hence by Corollary \ref{fcomp} we have $G^{\hat\ba}((\varnothing,\varnothing,(1)))=\ket{(\varnothing,\varnothing,(1))}$.

Now the non-empty ladders of $(2,1)$ are $\lad_1,\lad_2,\lad_e$ of residues $\hat0,\hat{-\negthinspace1},\hat1$ respectively, each containing one node, so we compute
\begin{align*}
A&=f_{\hat1}f_{\hat{-\negthinspace1}}f_{\hat0}\ket{(\varnothing,\varnothing,(1))}\\
&= \ket{((2,1),\varnothing,(1))}+q\ket{((1^2),(1),(1))}+q^2\ket{((1^2),\varnothing,(2))}+q\ket{((2),\varnothing,(1^2))}+q^2\ket{((1),(1),(1^2))}+q^3\ket{((1),\varnothing,(2,1))},\\
\end{align*}
and we see that $G^{\hat{\ba}}(\mu)=A$, independently of the choice of $e$.
\end{eg}

\subs{The canonical basis for $M^\ba$}

The main interest in this paper is in computing the canonical basis for the irreducible highest-weight module $M^\ba$; we have computed the canonical basis for the larger module $M^{\otimes\ba}$ simply in order to allow a recursive construction to work.  In order to obtain the canonical basis for $M^\ba$, one simply discards the unneeded vectors.  Here we comment briefly on how to identify these vectors.

The canonical basis for $M^\ba$ consists of the canonical basis vectors labelled by a certain set of multipartitions called \emph{regular multipartitions} in \cite{bk} or \emph{conjugate Kleshchev multipartitions} in \cite{wt2}; this result follows from \cite[Corollary 2.11]{am}.  (The latter uses a different $\UU$-action on $\calf^\ba$ and a different tensor product on $\UU$-modules, and is therefore stated in terms of Kleshchev multipartitions, but the translation between the two conventions is straightforward.)  We do not define (conjugate) Kleshchev multipartitions here, because the definition can be found in several places; but we note that the definition is recursive (though there has been some recent progress \cite{akt} towards giving a non-recursive definition).

Therefore one can obtain the canonical basis for $M^\ba$ from that for $M^{\otimes\ba}$ by computing the list of regular multipartitions and discarding canonical basis vectors not labelled by these.  However, we conjecture that there is a way to do this without computing the list of regular multipartitions.  In \cite{weight}, the author defined the notion of the \emph{weight} of a multipartition; this is a non-negative integer which depends on the multipartition and on $\ba$ (and should not be confused with the Lie-theoretic notion of weight).  In \cite{wt2}, we then proved a theorem which shows how this weight function is manifested in canonical bases.  Specifically (writing $w(\mu)$ for the weight of $\mu$, and translating from Kleshchev to regular multipartitions), we have the following.

\begin{propn}\label{wtcb}\thmcite{wt2}{Corollary 2.4}
If $\mu$ is a regular multipartition, then there is a multipartition $\la$ such that $d^\ba_{\la\mu}=q^{w(\mu)}$, while $d^\ba_{\nu\mu}$ has degree less than $w(\mu)$ for any other multipartition $\nu$.
\end{propn}

We conjecture that a converse to this statement is true: namely that if $\mu$ is a multipartition which is not regular, then the degree of $d^\ba_{\nu\mu}$ is less than $w(\mu)$ for all $\nu$.  This statement is proved in the case $r=1$ in \cite[Proposition 3.7]{runner}.  If this conjecture is true, then it leads to faster way to compute the canonical basis for $M^\ba$: one computes the canonical basis for $M^{\otimes \ba}$, computes the weight of each multipartition (which is quicker in general than checking whether a multipartition is regular), and then discards those canonical basis vectors $G^\ba(\mu)$ in which all the coefficients have degree less than $w(\mu)$.  We note in passing that this would yield a new (though relatively slow) recursive definition of regular multipartitions.
\newpage

\appendix

\sect{Index of notation}\label{app}

Since there is a great deal of notation involved in this paper, we include an index here for the reader's convenience.

\vspace{\topsep}
\begin{tabular}{ll}
$\inv ab$&$\{a,a+1,\dots,b\}$\\
$I$&$\zez$\\
$\calp$&the set of all partitions\\
$\reg$&the set of all $e$-regular partitions\\
$\varnothing$&the partition $(0,0,\dots)$\\
$\dom$&the dominance order on multipartitions\\
$[\la]$&the Young diagram of $\la$\\
$\la^{\fkn}$&the multipartition obtained by adding the node $\fkn$ to $[\la]$\\
$\la_{\fkn}$&the multipartition obtained by removing the node $\fkn$ from $[\la]$\\
$\tru\mu$&the $(r-1)$-multipartition $(\mu^{(2)},\dots,\mu^{(r)})$, for $\mu\in\calp^r$\\
$\ste\nu$&the $r$-multipartition $(\varnothing,\nu^{(1)},\dots,\nu^{(r-1)})$, for $\nu\in\calp^{r-1}$\\
$\mu_0$&$\ste{(\tru\mu)}$\\
$\add_i(\la)$&the set of addable $i$-nodes of $\la$\\
$\rem_i(\la)$&the set of removable $i$-nodes of $\la$\\
$\lad_l$&$l$th ladder in $\bbn^2$\\
$\lad_l(\mu)$&$\lad_l\cap[\mu]$\\
$\UU$&the quantum group $U_q(\widehat{\mathfrak{sl}}_e)$\\
$e_i,f_i,q^h$&generators of $\UU$\\
$\La_i$ ($i\in I$)&fundamental weights\\
$V(\La)$&irreducible highest-weight $\UU$-module with highest weight $\La$\\
$\ba=(s_1,\dots,s_r)$&element of $I^r$\\
$\tru\ba$&$(s_2,\dots,s_r)$\\
$\tilde\ba=(\tilde s_1,\dots,\tilde s_r)$&element of $\bbz^r$ such that $s_k=\tilde s_k+e\bbz$ for each $k$\\
$\res_{\bbz}(i,j,k)$&integral residue of a node $(i,j,k)$ (depending on $\tilde\ba$)\\
$\calf^\ba$&the Fock space associated with $\ba\in I^r$\\
$M^\ba$&submodule of $\calf^\ba$ generated by $\ket{\varnothing^r}$\\
$M^{\otimes\ba}$&$M^{(s_1)}\otimes\dots\otimes M^{(s_r)}$\\
$\ket\la$&standard basis element of $\calf^\ba$\\
$b^\ba_{\la\mu}$&coefficient of $\ket\la$ in $\ol{\ket\mu}$\\
$G^\ba(\mu)$&canonical basis element\\
$\calf^{\tilde\ba}$&twisted Fock space associated with $\tilde\ba$\\
$d^\ba_{\la\mu}$&coefficient of $\ket\la$ in $G^\ba(\mu)$\\
$\kt{\la,\tilde\ba}$&ordered wedge corresponding to $\la$ and $\tilde\ba$\\
$\one$&set of integers whose residue modulo $er$ lies in $\inv1e$\\
$H(\mu)$&basis element for $M^{\otimes\ba}$\\
\end{tabular}

\end{document}